\documentclass{article}
\usepackage{amssymb}
\usepackage{amsthm}
\usepackage{amssymb}
\usepackage{amsthm}
\usepackage{hyperref}
\usepackage{mathrsfs}
\usepackage{amsfonts}
\usepackage{amsmath}
 \usepackage{graphicx}
 \usepackage{subfigure}

\newtheorem{remark}{Remark}[section]

\newtheorem{theorem}{Theorem}[section]
\newtheorem{Lemma}{Lemma}[section]
\newtheorem{proposition}{Proposition}[section]

\begin{document}
\title{Precise high moment asymptotics for parabolic Anderson model with log-correlated Gaussian field}
\author{Yangyang Lyu$^{a,}${\thanks{Correspondence to: School of Mathematics, Jilin University, Changchun 130012, China.}\footnote{E-mail address: lvyy16@mails.jlu.edu.cn(Yangyang Lyu).} }\\[8pt]
\small \emph{\it $a$ School of Mathematics, Jilin University, Changchun 130012, China}}
\date{}\maketitle

\maketitle{}\mbox{}
\begin{center}
{
\begin{minipage}{16cm}
{{\textbf{Abstract}\: In this paper, we consider the continuous parabolic Anderson model (PAM) driven by a time-independent log-correlated Gaussian field (LGF).
We obtain an asymptotic result of
 $$\mathbb{E}\exp\Bigg\{\frac{1}{2}\sum\limits_{ j,k=1}^N\int_0^t\int_0^t\gamma(B_j(s)-B_k(r))drds\Bigg\}\qquad(N\rightarrow
 \infty)$$
which is composed of the independent Brownian motions $\{B_j(s)\}$ and the function $\gamma$ approximating to a logarithmic potential at $0$, such as the covariances of massive free field and Bessel field. Based on the asymptotic result, we get the precise high moment asymptotics for Feynman-Kac formula of the PAM with LGF.\\[8pt]
{\small {\bf Keywords}}\: Precise high moment asymptotics; Large deviation; Log-correlated Gaussian field; Massive free field; Bessel field\\}}
\end{minipage}
}
\end{center}

\section{Introduction}

In this paper, we consider the continuous parabolic Anderson model (PAM)
\begin{eqnarray}
\qquad\qquad\quad\left\{
\begin{array}{l}
\frac{\partial u}{\partial t}(t,x)=\frac{1}{2}\triangle u(t,x)+u(t,x)V(x),\qquad (t,x)\in \mathrm{R}_+\times\mathrm{R}^{d},\\
u(0,x)=u_0(x),\qquad x\in \mathrm{R}^{d},
\end{array}\label{Para}
\right.
\end{eqnarray}
with the log-correlated Gaussian field $V(x)$ on $\mathrm{R}^{d}$,
which is given by a centered Gaussian family $\{\langle V,\varphi\rangle; \varphi\in\mathcal{S}(\mathrm{R}^d)\}$ with covariance
\begin{eqnarray}
\qquad\qquad\qquad\mathbb{E}[\langle V,\phi\rangle \langle V,\psi\rangle]=\int_{\mathrm{R}^d\times\mathrm{R}^d}\phi(x)\psi(y)k(x,y) dxdy \qquad\forall\phi, \psi\in\mathcal{S}(\mathrm{R}^d).\label{2019329142754}
\end{eqnarray}
Here, $ \mathcal{S}(\mathrm{R}^d)$ is Schwartz space, $k(x,y)$ is a positive definite kernel on $\mathrm{R}^{d}\times\mathrm{R}^{d}$ and
\begin{eqnarray}
k(x,y)=\log_+\frac{T}{|x-y|}+g(x,y),\label{20181018215909}
\end{eqnarray}
where $\log_+x=(\log x)\vee0$, $\log$ is the natural logarithm, the constant $T>0$ is called correlation length and $g(x,y)$ is a bounded function on $\mathrm{R}^{d}\times\mathrm{R}^{d}$. For the details of LGF, see e.g. \cite{S25,M25,R27,L31}.
Throughout the paper, we assume that the initial value $u_0(x)$ satisfies
\begin{eqnarray}
0<\inf\limits_{x\in\mathrm{R}^d}u_0(x)\le\sup\limits_{x\in\mathrm{R}^d}u_0(x)
<+\infty. \label{20181217092348}
\end{eqnarray}
Because Stratonovich integral is closer to ``physical integral'' than Skorokhod integral, we interpret the $uV$ in (\ref{Para}) in the Stratonovich's sense. For the details of Stratonovich and Skorokhod integral, see e.g. \cite{S14}.

There are many achievements for the continuous PAM with other random fields, such as Poisson potential, fractional noise, etc. (see e.g.\cite{S3,M7,Q20,F15,P16}).
In view of the fruitful developments of LGF in quantum field theory, Wolfgang K\"{o}nig proposed the study of the PAM with LGF in his book \cite{W18}. As far as we know, there are few results for the high moment asymptotics of solutions to the model.

The precise high moment asymptotic problem for the solution $u(t,x)$ is to find a suitable rate $\sigma(N)$ about positive integer $N$ such that the quantity
$
\frac{1}{\sigma(N)}\log\mathbb{E}u^N(t,x)
$
converges to a nonzero constant as $N\rightarrow\infty$.
To some extent, it is meaningful to study the problem in physics.
One hand, in \cite{M10,S11}, the precise high moment asymptotics is a key step in the process of computing spatial asymptotics. As explained in \cite{O9}, spatial asymptotics is one of the problems of ``physical intermittency"
 which is applied to describing the absence of diffusion of waves in a disordered medium in physics.
 On the other hand, based on the precise high moment asymptotic result, the precise asymptotics of tail probability of $\log u(t,x)$ can be obtained by nonnegative large deviation technique in \cite{S19}. Here, the $\log u(t,x)$ represents the free energy of directed Brownian polymer in a continuous random environment.
 The high moment asymptotic problem has attracted a lot of attention.
There have been many results for precise high moment asymptotics
when the covariances of driven Gaussian fields are homogeneous in time and space, or locally integrable in time and bounded in space, see \cite{F1,M10,S11,P16,S19,P21} as references.

 By the following (\ref{201943095707}) and the same method as \cite{S14}, it's trivial to prove that if the Feynman-Kac formula of equation (\ref{Para}) exists then it is a mild solution in the Stratonovich's sense. However, according to \cite{S14,F15}, the uniqueness of the mild and weak solutions of PAM can't be proved in the Stratonovich's sense now.
 Thus, this paper is only concerned with the precise high moment asymptotics of the Feynman-Kac formula of equation (\ref{Para}). We remark that,
 in the Young's sense, it can be proved that the Feynman-Kac formula is the unique mild solution for PAM with the Gaussian fields in \cite{S14,Q35}. The similar result might be also extended to the PAM with LGF, which will be our future work.

According to \cite{F15}, the Feynman-Kac formula for the solution to equation (\ref{Para}) is
\begin{eqnarray}
u(t,x)=\mathbb{E}^{B}\left[\exp\left\{\int_0^tV(B^x(s))ds\right\}u_0(B^x(t))
\right],\label{2019324195809}
\end{eqnarray}
where $B^x$ is a $d$-dimensional Brownian motion independent of $V$ and  starting at $x\in \mathrm{R}^{d}$ and $\mathbb{E}^{B}$ is the expectation with respect to the standard Brownian motion $B$.
For all $\varepsilon>0$, let $p_\varepsilon(x):=(2\pi\varepsilon)^{-\frac{d}{2}}\exp\{-\frac{1}
{2\varepsilon}|x|^2\}$ and $V_\varepsilon(x):=\langle V(\cdot), p_\varepsilon(x-\cdot)\rangle$, then for all $t>0$, the time integral in (\ref{2019324195809}) is defined as the $L^2$-limit
\begin{eqnarray}
\int_0^tV(B(s))ds:=\lim\limits_{\varepsilon\downarrow0}\int_0^t
V_{\varepsilon}(B(s))ds-\mbox{ $L^2(\Omega,\mathcal{F},\mathbb{P})$}.\label{2019323173155}
\end{eqnarray}
It's easy to get that (\ref{2019323173155}) is well defined by (\ref{20181018215909}).
Conditioning on the Brownian motion, the $\int_0^tV(B(s))ds$ is a centered Gaussian process with the conditional covariance $\int_0^t\int_0^tk(B(s),B(r))dsdr$. Furthermore, by reference to \cite{F15},  the Feynman-Kac moment representation of (\ref{2019324195809}) is
\begin{eqnarray}
\mathbb{E}u^N(t,x)
=\mathbb{E}\Bigg[\exp\Bigg\{\frac{1}{2}\sum\limits_{ j,k=1}^N\int_0^t\int_0^tk(B_j^x(s),B_k^x(r))drds\Bigg\}\prod\limits_{j=1}^Nu_0(B_j^x(t))
\Bigg],\label{20181226223826}
\end{eqnarray}
where $\{B_{j}^x; j=1,\cdots,N\}$ is a family of $d$-dimensional independent B.M. starting at $x$. By reference to \cite{F15}, the existence of (\ref{2019324195809}) and (\ref{20181226223826}) is provided by the exponential integrability of $\int_0^t\int_0^tk(B(s),B(r))dsdr$.
In fact, by (\ref{20181018215909}), there exists some $\alpha\in(0,1)$ and $C>0$ such that
\begin{eqnarray}
|k(x,y)|\le|x-y|^{-\alpha}+C,\label{201943095707}\qquad\qquad\qquad\qquad\qquad
\qquad\qquad\qquad\\
\mbox{then\qquad\qquad}\mathbb{E}\exp
\left\{\int_0^t\int_0^tk(B(s),B(r))dsdr\right\}\le e^{Ct^2}\mathbb{E}\exp
\left\{\int_0^t\int_0^t|B(s)-B(r)|^{-\alpha}dsdr\right\}.\label{2019330183219}
\qquad\qquad
\end{eqnarray}
Hence, the exponential integrability of $\int_0^t\int_0^tk(B(s),B(r))dsdr$ is transformed into that of the right-hand side in (\ref{2019330183219}).
To explain the exponential integrability of $\int_0^t\int_0^t|B(s)-B(r)|^{-\alpha}dsdr$ and prove the following Proposition \ref{C1}, we need to introduce a general framework.
Let $\gamma$ be the Fourier transform of the tempered measure $\mu$ in $\mathcal{S}'(\mathrm{R}^d)$, i.e.
\begin{eqnarray}
\qquad\qquad\qquad\langle\gamma, \phi\rangle=\int_{\mathrm{R}^d}\mathcal{F}(\phi)(\xi)\mu(d\xi),\qquad\qquad\forall\phi\in
\mathcal{S}(\mathrm{R}^d).\label{tmm}
\end{eqnarray}
Here, $\mathcal{F}(\phi)(\xi):=\int_{\mathrm{R}^d}e^{i\xi\cdot x}\phi(x)dx $ is the Fourier transform of $\phi$ on $\mathrm{R}^d$.
For all $\varepsilon>0$, let $\mu_\varepsilon(d\xi):=\exp\left\{-\frac{\varepsilon}{2}|\xi|^2\right\} \mu(d\xi)$ and $\gamma_\varepsilon(x):=\int_{\mathrm{R}^d}e^{i\xi\cdot x}\mu_\varepsilon(d\xi)$, then we define
\begin{eqnarray}\label{5gun}
\int_0^t\int_0^t\gamma(B(s)-B(r))dsdr
:=\lim\limits_{\varepsilon\downarrow0}\int_0^t
\int_0^t\gamma_{\varepsilon}(B(s)-B(r))dsdr-\mbox{ $L^1(\Omega,\mathcal{F},\mathbb{P})$}.\label{2019314224354}
\end{eqnarray}
By (\ref{2019314224354}) and Proposition 4.4 in \cite{Q35}, it can be obtained that the Dalang's condition
$
\int_{\mathrm{R}^d}\frac{1}{1+|\xi|^2}\mu(d\xi)<+\infty
$ is equivalent to
\begin{eqnarray}
\mathbb{E}\exp
\left\{\int_0^t\int_0^t\gamma(B(s)-B(r))dsdr\right\}<+\infty.
\label{11111}
\end{eqnarray}
Becasue the Riesz potential $|\cdot|^{-\alpha}$ in (\ref{2019330183219}) satisfies the Dalang's condition, $\int_0^t\int_0^tk(B(s),B(r))dsdr$ is exponentially integrable. Hence, (\ref{2019324195809}) and (\ref{20181226223826}) are well defined. Based on it, we obtain
the following main result about the precise high moment asymptotics for the Feynman-Kac formula (\ref{2019324195809}) of the PAM with LGF.
\begin{theorem}\label{201946171643}
For the Feynman-Kac formula (\ref{2019324195809}), it has
\begin{eqnarray}
\lim\limits_{N\rightarrow\infty}\frac{1}{N \lambda_N}\log\mathbb{E}u^N(t,x)
=\frac{t^2}{2}.\label{201942104514}
\end{eqnarray}
Here, $\lambda_N$ is a function about positive integer $N$, which satisfies $\lambda_N>e$ and $\frac{\log \lambda_N}{\lambda_N} =\frac{2}{N}$ when $N$ is enough large.
\end{theorem}
\begin{remark}\label{2019412145239}
There is usually no explicit solution for the equation $\frac{\log \lambda_N}{\lambda_N} =\frac{2}{N}$.
When $N$ is enough large, the equation has two implicit solutions, and the $\lambda_N$ is unique if we assume $\lambda_N>e$. Moreover, the $\lambda_N$ is a regularly varying function about $N$ and satisfies $\frac{1}{\lambda_N}=o(\frac{1}{N})$ as $N\rightarrow\infty$.
\end{remark}
The paper is organized as follows. In Section \ref{20181212231613}, we prove Proposition \ref{C1} which is applied to proving Proposition \ref{2018109231801}. In Section \ref{20181225200911}, we prove Proposition \ref{2018109231801} which is an asymptotic result for the function approximating to a logarithmic potential at $0$. Based on Proposition \ref{2018109231801}, Section \ref{201947120052} is the proof of Theorem \ref{201946171643}.
\section{A general limit in the time-independent case}\label{20181212231613}

To prove Proposition \ref{2018109231801}, we need the following Proposition \ref{C1} which is a general result for positive definite function in the time-independent case. In view of the time-space case, i.e. $\sum_{j,k=1}^N\int_0^{t}\int_0^{t}\gamma_0(s-r)\gamma(B_j(s)-B_k(r))dsdr$, we call $\sum_{j,k=1}^N\int_0^{t}\int_0^{t}\gamma(B_j(s)-B_k(r))dsdr$ and $\sum_{1\le j\neq k\le N}\int_0^{t}\gamma(B_j(s)-B_k(s))ds$ time-independent case ($\gamma_0\equiv1$) and time-white case ($\gamma_0=\delta_0$, i.e. Dirac delta function at $0$), respectively. Before Proposition \ref{C1}, the time-white case has been considered in \cite{F1} . When the $\mu$ is finite in Proposition \ref{C1}, our method is partly inspired by \cite{F1} and we use the semigroup inequality in \cite{R2}. When the $\mu$ is unbounded, we can't transform the time-independent case into the time-white case because the $\sum_{1\le j\neq k\le N}\int_0^{t}\gamma(B_j(s)-B_k(s))ds$ may be not exponentially integrable.
\begin{proposition}\label{C1}
Assume that the $\gamma$ in (\ref{tmm}) satisfies the Dalang's condition, then for any sequence $\{t_N\}_{N\in\mathbb{N}}$ satisfying $t_N\rightarrow+\infty$ as $N\rightarrow\infty$, it holds that
\begin{eqnarray}
&&\lim\limits_{N\rightarrow\infty}\frac{1}{Nt_N}\log\mathbb{E}\exp
\Bigg\{\frac{1}{2Nt_N}\sum\limits_{j, k=1}^N\int_0^{t_N}\int_0^{t_N}\gamma(B_j(s)-B_k(r))dsdr\Bigg\}\label{201946101507}\\
&&=\frac{1}{2}\sup\limits_{g\in\mathcal{F}_d}
\Bigg\{\int_{\mathrm{R}^d}\Bigg|\int_{\mathrm{R}^d}e^{i\xi\cdot x}g^2(x)dx\Bigg|^2\mu(d\xi)
-\int_{\mathrm{R}^d}|\nabla g(x)|^2dx\Bigg\}:=\mathcal{E}(\gamma)\label{X2}
\end{eqnarray}
and the variation $\mathcal{E}(\gamma)<+\infty$,  where we define the set of functions
$
\mathcal{F}_d:=\left\{g\in W^{1,2}(\mathrm{R}^d);\|g\|_2=1\right\}
$.
\end{proposition}

To prove Proposition \ref{C1}, we need the following two Lemmas.
\begin{Lemma}\label{C3.1}
For any $t_1,t_2>0$ and positive integer $m$, it holds that
\begin{eqnarray}
~\mathbb{E}\exp\Bigg\{\frac{1}{t_1+ t_2}\sum_{j,k=1}^m\int_0^{ t_1 + t_2}\int_0^{ t_1 + t_2}
\gamma(B_j(s)-B_k(r))dsdr\Bigg\}\qquad\qquad\qquad\qquad\qquad\qquad\qquad
\qquad\qquad\nonumber\\
~\le\mathbb{E}\exp\Bigg\{\frac{1}{ t_1 }\sum_{j,k=1}^m\int_0^{ t_1 }\int_0^{ t_1 }
\gamma(B_j(s)-B_k(r))dsdr\Bigg\}\cdot\mathbb{E}\exp\Bigg\{\frac{1}{ t_2}\sum_{j,k=1}^m\int_0^{ t_2}\int_0^{ t_2}
\gamma(B_j(s)-B_k(r))dsdr\Bigg\}.\label{2019315102011}
\end{eqnarray}
\end{Lemma}

\noindent\textbf{proof:}
By the approximation in (\ref{2019314224354}) and Jensen's inequality, the following inequality is almost surely established.
\begin{eqnarray}
&&\frac{1}{t_1+t_2 }\int_{\mathrm{R}^{d}}\Bigg|\sum_{j=1}^m\int_{0}^{t_1+t_2}e^{i\xi\cdot B_j(s)}ds\Bigg|^2\mu(d\xi)\nonumber\\
&&\le \frac{1}{t_1}\int_{\mathrm{R}^{d}}\Bigg|\sum_{j=1}^m\int_{0}^{t_1}e^{i\xi\cdot B_j(s)}ds\Bigg|^2\mu(d\xi)+\frac{1}{t_2}\int_{\mathrm{R}^{d}}\Bigg|\sum_{j=1}^m
\int_{t_1}^{t_1+t_2 }
e^{i\xi\cdot B_j(s)}ds\Bigg|^2\mu(d\xi).\label{201315101645}
\end{eqnarray}
Moreover, by (\ref{201315101645}) and strong Markov property, we have
\begin{eqnarray}
&&\mathbb{E}\exp\Bigg\{\frac{1}{t_1+ t_2}\sum_{j,k=1}^m\int_0^{ t_1 + t_2}\int_0^{ t_1 + t_2}
\gamma(B_j(s)-B_k(r))dsdr\Bigg\}\nonumber\\
&&\le\mathbb{E}\Bigg[\exp\Bigg\{\frac{1}{ t_1 }\sum_{j,k=1}^m\int_0^{ t_1 }\int_0^{ t_1 }
\gamma(B_j(s)-B_k(r))dsdr\Bigg\}\nonumber\\
&&\cdot\mathbb{E}\Bigg[\exp\Bigg\{\frac{1}{ t_2}\sum_{j,k=1}^m\int_{t_1 }^{t_1 + t_2}\int_{t_1 }^{ t_1 +t_2}
\gamma(B_j(s)-B_k(r))dsdr\Bigg\}\Bigg|\mathfrak{F}_{t_1}\Bigg]\Bigg]
\nonumber\\
&&\le\mathbb{E}\exp\Bigg\{\frac{1}{ t_1 }\sum_{j,k=1}^m\int_0^{ t_1 }\int_0^{ t_1 }
\gamma(B_j(s)-B_k(r))dsdr\Bigg\}\nonumber\\
&&\cdot\Bigg(\sup\limits_{\tilde{b}}\mathbb{E}^{\tilde{b}}\exp\Bigg\{\frac{1}{ t_2}\sum_{j,k=1}^m
\int_{0 }^{ t_2}\int_{0}^{ t_2}
\gamma(B_j(s)-B_k(r))dsdr\Bigg\}\Bigg),\label{2019315121809}
\end{eqnarray}
where $\tilde{b}=(b_1,\cdots,b_m)\in (R^{d})^m$ and $\mathfrak{F}_{t_1}$ is the $\sigma$-algebra up to $t_1$.
We claim that for any $t>0$, it holds that
\begin{eqnarray}
\sup\limits_{\tilde{b}}\mathbb{E}^{\tilde{b}}\exp\Bigg\{\frac{1}{ t}\sum_{j,k=1}^m
\int_{0 }^{ t}\int_{0}^{ t}
\gamma(B_j(s)-B_k(r))dsdr\Bigg\}\le\mathbb{E}\exp\Bigg\{\frac{1}{ t}\sum_{j,k=1}^m
\int_{0 }^{ t}\int_{0}^{ t}
\gamma(B_j(s)-B_k(r))dsdr\Bigg\}.\label{2019315153606}
\end{eqnarray}
Indeed, by Taylor expansion and Bochner's representation, for all integer $n\ge1$ and $\tilde{b}\in (R^{d})^m$, we have
\begin{eqnarray}
&&\mathbb{E}^{\tilde{b}}\Bigg[\frac{1}{ t}\sum_{j,k=1}^m
\int_{0 }^{ t}\int_{0}^{ t}
\gamma(B_j(s)-B_k(r))dsdr\Bigg]^n\nonumber\\
&&=\frac{1}{t^n}\int_{\mathrm{R}^{dn}}\int_{[0,t]^{2n}}
\sum_{ j_1,\cdots,j_n=1\atop k_1,\cdots,k_n=1
}^m\mathbb{E}\prod\limits_{l=1}^n\mathrm{e}^{i\xi_l \cdot (b_{k_l}-b_{j_l})}\cdot\mathrm{e}^{i\xi_l\cdot (B_{k_l}(s_l)-B_{j_l}(r_l))}
\prod\limits_{l=1}^nds_l
dr_l
\mu(d\xi_l)\nonumber\\
&&\le\frac{1}{t^n}\int_{\mathrm{R}^{dn}}\int_{[0,t]^{2n}}
\sum_{ j_1,\cdots,j_n=1\atop k_1,\cdots,k_n=1
}^m\mathbb{E}\mathrm{e}^{i\sum\limits_{l=1}^n\xi_l\cdot (B_{k_l}(s_l)-B_{j_l}(r_l))}
\prod\limits_{l=1}^nds_l
dr_l
\mu(d\xi_l)\nonumber\\
&&=\mathbb{E}\Bigg[\frac{1}{ t}\sum_{j,k=1}^m
\int_{0 }^{ t}\int_{0}^{ t}
\gamma(B_j(s)-B_k(r))dsdr\Bigg]^n.\label{2019315155230}
\end{eqnarray}
Then, by (\ref{2019315121809}) and (\ref{2019315153606}), the proof is completed.

The following Lemma is a direct result of Lemma \ref{C3.1}.
\begin{Lemma}\label{C3.2}
For any positive integer $m$, $n$ and $t>0$, it holds that
\begin{eqnarray}
\mathbb{E}\exp
\Bigg\{\frac{1}{n}\sum_{j,k=1}^m\int_0^{nt}\int_0^{nt}
\gamma(B_j(s)-B_k(r))dsdr\Bigg\}\le\Bigg(\mathbb{E}\exp
\Bigg\{\sum_{j,k=1}^m\int_0^{t}\int_0^{t}
\gamma(B_j(s)-B_k(r))dsdr\Bigg\}\Bigg)^n.\label{2019315164245}
\end{eqnarray}
\end{Lemma}

\noindent\textbf{The proof of Proposition \ref{C1}:}
 \textit{Step 1.} We prove (\ref{X2}) in the case that $\mu$ is finite.
Firstly, we prove the lower bound of (\ref{X2}). When $\mu$ is finite, by Bochner's representation, we have
\begin{eqnarray}
&&\mathbb{E}\exp
\Bigg\{\frac{1}{2Nt_N}\sum\limits_{j, k=1}^N\int_0^{t_N}\int_0^{t_N}\gamma(B_j(s)-B_k(r))dsdr\Bigg\}\nonumber\\
&&=\mathbb{E}\exp\Bigg\{\frac{Nt_N}{2}\int_{\mathrm{R}^d}\Bigg|
\frac{1}{Nt_N}\sum\limits_{k=1}^N
\int_0^{t_N}\mathrm{e}^{i\xi \cdot B_k(s)}ds\Bigg|^2\mu(d\xi)\Bigg\}.
\label{2018119154309}\label{2018119154327}
\end{eqnarray}
  Let $\mathcal{H}$ be a subspace of complex-value Hilbert space $L^2(\mathrm{R}^d,\mu)$ satisfying $f(-x)=\overline{f(x)}$ and $\langle f,h\rangle_{\mathcal{H}}:=\int_{\mathrm{R}^d}f(\xi)\overline{h(\xi)}
 \mu(d\xi)$ for all $ f,h\in\mathcal{H}$. Then, by (\ref{2018119154309}), $\|h\|_{\mathcal{H}}^2\ge 2\langle f,h\rangle_{\mathcal{H}}-\|f\|_{\mathcal{H}}^2$ for all $f,h\in\mathcal{H}$ and Theorem 4.1.6 in \cite{R2}, we have
\begin{eqnarray}
&&\liminf\limits_{N\rightarrow\infty}\frac{1}{Nt_N}\log\mathbb{E}\exp
\Bigg\{\frac{1}{2Nt_N}\sum\limits_{j, k=1}^N\int_0^{t_N}\int_0^{t_N}\gamma(B_j(s)-B_k(r))dsdr\Bigg\}\nonumber\\
&&\ge\sup\limits_{f\in\mathcal{H}}\left\{\liminf\limits_{N\rightarrow\infty}\frac{1}{t_N}\log\mathbb{E}
\exp\left\{\int_0^{t_N}f^*(B(s))ds\right\}-\frac{1}{2}\|f\|_{\mathcal{H}}^2\right\},\label{4}\\
&&=\sup\limits_{f\in\mathcal{H}}\sup\limits_{g\in\mathcal{F}_d}
\left\{\int_{\mathrm{R}^d}f^*(x)g^2(x)dx
-\frac{1}{2}\int_{\mathrm{R}^d}|\nabla g(x)|^2dx-\frac{1}{2}\|f\|_{\mathcal{H}}^2\right\}\nonumber\\
&&=\sup\limits_{g\in\mathcal{F}_d}\sup\limits_{f\in\mathcal{H}}
\left\{\langle f,\mathcal{F}(g^2)\rangle_{\mathcal{H}}-\frac{1}{2}\|f\|_{\mathcal{H}}^2
-\frac{1}{2}\int_{\mathrm{R}^d}|\nabla g(x)|^2dx\right\}=\mathcal{E}(\gamma),\label{20181223223240}
\end{eqnarray}
where $f^*(x)=\int_{\mathrm{R}^d}e^{i\xi\cdot x}f(-\xi)\mu(d\xi)$ and the last equality is due to $\|h\|_{\mathcal{H}}^2=\sup\limits_{f\in\mathcal{H}}\{ 2\langle f,h\rangle_{\mathcal{H}}-\|f\|_{\mathcal{H}}^2\}$.

Secondly, we prove the upper bound of (\ref{X2}) when $\mu$ is finite. We first assume that $t_N/t$ is always an integer for all positive integer $N$ and fixed $t>0$.
By Lemma \ref{C3.2} and Bochner's representation, it holds that
\begin{eqnarray}
&&\frac{1}{Nt_N}\log\mathbb{E}\exp
\Bigg\{\frac{1}{2Nt_N}\sum\limits_{j, k=1}^N\int_0^{t_N}\int_0^{t_N}\gamma(B_j(s)-B_k(r))dsdr\Bigg\}\nonumber\\
&&\le\frac{1}{Nt}\log\mathbb{E}\exp\Bigg\{\frac{Nt}{2}\int_{\mathrm{R}^d}\left|\frac{1}{N}\sum\limits_{k=1}^N\frac{1}{t}\int_0^{t}\mathrm{e}^{i\xi \cdot B_k(s)}ds\right|^2\mu(d\xi)\Bigg\}.\label{2018224102701}
\end{eqnarray}
For all positive integer $1\le k\le N$, define $X_k:=\frac{1}{t}\int_0^{t}\mathrm{e}^{i\xi \cdot B_k(s)}ds$  and $X\stackrel{d}{=}X_k$. Let $\left\{X, X_k\right\}_{1\le k\le N}$ be a sequence of i.i.d $\mathcal{H}$-valued random variables. Apparently, the $X$ satisfies the exponential moment condition:
$
\mathbb{E}\exp\left\{\theta\left\|X\right\|_{\mathcal{H}}\right\}<\infty$ for all $ \theta>0.
$
Let Cram\'{e}r function and its Legendre transform be
$$\Lambda(y):=\log\mathbb{E}\exp\left\{\left\langle y,X\right\rangle_{\mathcal{H}}\right\} \quad \forall y\in\mathcal{H}\quad \mbox{and}\quad\Lambda^*(x):=\sup\limits_{y\in\mathcal{H}}\left\{\langle y,x\rangle_{\mathcal{H}}-\Lambda(y)\right\}\quad \forall x\in\mathcal{H},$$
respectively.
By Donsker-Varadhan's large deviation principle (LDP) in a separable Banach space (Theorem 5.3,\cite{D5}), the distributions  $\Big\{\mathbb{P}\Big(\frac{1}{N}\sum\limits_{k=1}^NX_k\in\cdot\Big)\Big\}_{N\in
\mathbb{N}}$ satisfy LDP in $\mathcal{H}$ with the rate function $\Lambda^*(x)$. By $\|X_k\|_{\mathcal{H}}^2\le \mu(\mathrm{R}^d)$ for all $ 1\le k\le N$ and Varadhan integral Lemma (theorem 4.3.1,\cite{D6}), we have
\begin{eqnarray}
\lim\limits_{N\rightarrow\infty}\frac{1}{N}\log\mathbb{E}\exp\Bigg\{\frac{Nt}{2}
\int_{\mathrm{R}^d}\Bigg|\frac{1}{N}\sum\limits_{k=1}^N
\frac{1}{t}\int_0^{t}\mathrm{e}^{i\xi \cdot B_k(s)}ds\Bigg|^2\mu(d\xi)\Bigg\}
=\sup\limits_{h\in\mathcal{H}}\left\{\frac{t}{2}\|h\|_{\mathcal{H}}^2
-\Lambda^*(h)\right\}.\label{26}
\end{eqnarray}
Here, we claim that for all $t>0$, it holds that
\begin{eqnarray}
\sup\limits_{h\in\mathcal{H}}\left\{\frac{t}{2}\|h\|_{\mathcal{H}}^2-\Lambda^*(h)\right\}
=\sup\limits_{\|h\|_{\mathcal{H}}^2\le\mu(\mathrm{R}^d)}\left\{\frac{t}{2}\|h\|_{\mathcal{H}}^2-\Lambda^*(h)\right\}.\label{20181030160814}
\end{eqnarray}
In fact, the
 $\frac{1}{N}\sum\limits_{k=1}^NX_k$ satisfies lower bound of LDP in $\mathcal{H}$. For the open set $[\|\cdot\|_{\mathcal{H}}^2>\mu(\mathrm{R}^d)]$ of $\mathcal{H}$, we have
\begin{eqnarray}
\liminf\limits_{N\rightarrow\infty}\frac{1}{N}\log\mathbb{P}
\Bigg(\Bigg\|\frac{1}{N}
\sum\limits_{k=1}^NX_k\Bigg\|_{\mathcal{H}}^2>\mu(\mathrm{R}^d)\Bigg)\ge
-\inf\limits_{\|h\|_{\mathcal{H}}^2>\mu(\mathrm{R}^d)}\Lambda^*(h).\label{20181226164738}
\end{eqnarray}
By (\ref{20181226164738}) and the fact $\Big\|\frac{1}{N}\sum\limits_{k=1}^NX_k\Big\|_{\mathcal{H}}^2\le\mu(\mathrm{R}^d)$,
we get $\Lambda^*(h)=+\infty$ on the set $\{h\in\mathcal{H};\|h\|_{\mathcal{H}}^2>\mu(\mathrm{R}^d)\}$. Hence, the quantity in (\ref{20181030160814}) reaches the maximum on the set $\{h\in\mathcal{H};\|h\|_{\mathcal{H}}^2\le\mu(\mathrm{R}^d)\}$.
We notice that, for all $h\in\mathcal{H}$ and $\|h\|_{\mathcal{H}}^2\le\mu(\mathrm{R}^d)$, the linearization of $h$ can be rewritten as
\begin{eqnarray}
\|h\|_{\mathcal{H}}^2=\sup\limits_{\|f\|_{\mathcal{H}}^2\le\mu(\mathrm{R}^d) }\{ 2\langle f,h\rangle_{\mathcal{H}}-\|f\|_{\mathcal{H}}^2\},\label{20181226172706}
\end{eqnarray}
which is due to that the quantity in (\ref{20181226172706}) attains the maximum when $f=h$.
By (\ref{20181030160814}), (\ref{20181226172706}) and Legendre transform, we get
\begin{eqnarray}
\sup\limits_{\|h\|_{\mathcal{H}}^2\le\mu(\mathrm{R}^d)}\left\{\frac{t}{2}\|h\|_{\mathcal{H}}^2-\Lambda^*(h)\right\}
&=&\sup\limits_{\|h\|_{\mathcal{H}}^2\le\mu(\mathrm{R}^d)}
\Bigg\{\sup\limits_{\|f\|_{\mathcal{H}}^2\le\mu(\mathrm{R}^d)}
\left\{-\frac{t}{2}\|f\|_{\mathcal{H}}^2+t\langle f,h\rangle_{\mathcal{H}}\right\}-\Lambda^*(h)\Bigg\}\nonumber\\
&\le&
\sup\limits_{\|f\|_{\mathcal{H}}^2\le\mu(\mathrm{R}^d)}\left\{-\frac{t}{2}\|f\|_{\mathcal{H}}^2+\Lambda(tf)\right\}.\label{20181224102805}
\end{eqnarray}
Let $I(x):=\langle f,\mathrm{e}^{i\xi\cdot x}\rangle_{\mathcal{H}}$ for all $ f\in\mathcal{H}$, then $I(x)$ is a bounded and continuous function about $x$ on $\mathrm{R}^d$. By (\ref{2018224102701}), (\ref{26}) and (\ref{20181224102805}), we only need to prove
\begin{eqnarray}\label{30}
&&\limsup\limits_{t\rightarrow+\infty}\sup\limits_{\|f\|_{\mathcal{H}}^2
\le\mu(\mathrm{R}^d)}\frac{1}{t}
\log\mathbb{E}\exp\left\{-\frac{t}{2}\|f\|_{\mathcal{H}}^2+\int_0^tI(B(s))ds\right\} \label{201941170040}\\
&&\leq\frac{1}{2}\sup\limits_{g\in\mathcal{F}_d}\left\{\int_{\mathrm{R}^d\times\mathrm{R}^d}\gamma(x-y)g^2(x)g^2(y)dxdy-\int_{\mathrm{R}^d}|\nabla g(x) |^2dx\right\}.\nonumber
\end{eqnarray}
Indeed,
let $C_1:=\sup\limits_{x\in\mathrm{R}^{d} \atop \|f\|_{\mathcal{H}}^2\le\mu(\mathrm{R}^d) }
|I(x)|\le\mu(\mathrm{R}^d)$, then there exists positive constants $C_2, C_3, C_4<+\infty$ such that for all $t>1$, it has
\begin{eqnarray}
&&\mathbb{E}\exp\left\{\int_0^tI(B(s))ds\right\}\le e^{C_1}\mathbb{E}\exp\left\{\int_1^tI(B(s))ds\right\}\le e^{C_1}\Bigg(\mathbb{E}\Bigg[1_{|B(1)|\le t^2}\exp\left\{\int_1^tI(B(s))ds\right\}1_{|B(t)|\le t^2}\Bigg]\nonumber\\
&&+\mathbb{E}\Bigg[\exp\left\{\int_1^tI(B(s))ds\right\}1_{|B(1)|>t^2}\Bigg]
+\mathbb{E}\Bigg[\exp\left\{\int_1^tI(B(s))ds\right\}1_{|B(t)|> t^2}\Bigg]\Bigg)\nonumber\\
&&\le C_2e^{C_1}t^{2d}\exp\Bigg\{(t-1)\sup\limits_{g\in \mathcal{F}_d}\Big\{\int_{\mathrm{R}^d}I(x)g^2(x)dx-\frac{1}{2}\int_{\mathrm{R}^d}|\nabla g(x) |^2dx\Big\}\Bigg\}+ C_3e^{C_1t-\frac{t^4}{2}}+ C_4e^{C_1t-\frac{t^3}{2}},\label{20181226185624}
\end{eqnarray}
which has been proved in page 101-102 of \cite{R2}. Now we use the semigroup inequality (\ref{20181226185624}), then
\begin{eqnarray}
&&\sup\limits_{\|f\|_{\mathcal{H}}^2\le\mu(\mathrm{R}^d)}\mathbb{E}\exp\left\{-\frac{t}{2}\|f\|_{\mathcal{H}}^2+
\int_0^tI(B(s))ds\right\}\nonumber\\
&&\le C_2e^{C_1}t^{2d}\exp\Bigg\{(t-1)\sup\limits_{g\in \mathcal{F}_d}\Big\{\sup\limits_{  f\in\mathcal{H}}\Big\{-\frac{1}{2}\|f\|_{\mathcal{H}}^2+\Big\langle f(\xi),\int_{\mathrm{R}^d}\mathrm{e}^{i\xi\cdot x}g^2(x)dx\Big\rangle_{\mathcal{H}}
\Big\}-\frac{1}{2}\int_{\mathrm{R}^d}|\nabla g(x) |^2dx\Big\}\Bigg\}\nonumber\\
&&+C_3e^{C_1t-\frac{t^4}{2}}
+ C_4e^{C_1t-\frac{t^3}{2}}\nonumber\\
&&=C_2e^{C_1}t^{2d}\exp\left\{\frac{(t-1)}{2}\sup\limits_{g\in\mathcal{F}_d}
\left\{\int_{\mathrm{R}^d\times\mathrm{R}^d}\gamma(x-y)g^2(x)g^2(y)dxdy-
\int_{\mathrm{R}^d}|\nabla g(x) |^2dx\right\}\right\}+C_3e^{C_1t-\frac{t^4}{2}}\nonumber\\
&&
+ C_4e^{C_1t-\frac{t^3}{2}}.\label{20181226201635}
\end{eqnarray}
Then, by (\ref{20181226201635}) and the inequality $\log(a+b+c)\le\log3+\max\{\log a,\log b,\log c\}$ for all $a,b,c>0$, we have
\begin{eqnarray}
&&\limsup\limits_{t\rightarrow+\infty}\sup
\limits_{\|f\|_{\mathcal{H}}^2\le\mu(\mathrm{R}^d)}
\frac{1}{t}\log\mathbb{E}\exp\left\{-\frac{t}{2}
\|f\|_{\mathcal{H}}^2+\int_0^tI(B(s))ds\right\} \nonumber\\
&&\le\limsup\limits_{t\rightarrow+\infty}\max
\Bigg\{\frac{1}{2}\sup\limits_{g\in\mathcal{F}_d}\left\{\int_{\mathrm{R}^d\times\mathrm{R}^d}\gamma(x-y)g^2(x)g^2(y)dxdy-\int_{\mathrm{R}^d}|\nabla g(x) |^2dx\right\}+\frac{\log (C_2e^{C_1}t^{2d})}{t}, \nonumber\\
&&C_1-\frac{t^3}{2}+\frac{\log C_3}{t}, C_1-\frac{t^2}{2}+\frac{\log C_4}{t}\Bigg\}=\mathcal{E}(\gamma).\label{20181212201340}
\end{eqnarray}
We prove the result for the general $t_N/t$. For any parameters $ p,q>0$ and $ p+q=1$,
by Jensen's inequality, triangle inequality, (\ref{2018224102701}) and (\ref{20181212201340}), we have
\begin{eqnarray}
&&\limsup\limits_{N\rightarrow\infty}\frac{1}{Nt_N}\log\mathbb{E}\exp
\Bigg\{\frac{1}{2Nt_N}\sum\limits_{j, k=1}^N\int_0^{t_N}\int_0^{t_N}\gamma(B_j(s)-B_k(r))dsdr\Bigg\}\nonumber\\
&&\le\limsup\limits_{p\rightarrow1}\limsup\limits_{t\rightarrow+\infty}
\limsup\limits_{N\rightarrow\infty}
\frac{1}{N
\lfloor\frac{t_N}{t}\rfloor t}\log
\mathbb{E}\exp\Bigg\{\frac{1}{2pN\lfloor\frac{t_N}{t}\rfloor t}\int_{\mathrm{R}^d}
\left|\sum\limits_{k=1}^N\int_0^{\lfloor\frac{t_N}{t}\rfloor t}
\mathrm{e}^{i\xi \cdot B_k(s)}ds\right|^2\mu(d\xi)+\frac{Nt\mu(\mathrm{R}^d)}{2qt_N}\Bigg\}\nonumber\\
&&=\lim\limits_{p\rightarrow1}
\frac{1}{2}\sup\limits_{g\in\mathcal{F}_d}\left\{\int_{\mathrm{R}^d\times
\mathrm{R}^d}
\frac{1}{p}\gamma(x-y)g^2(x)g^2(y)dxdy-\int_{\mathrm{R}^d}|\nabla g(x) |^2dx\right\}=\mathcal{E}(\gamma).\label{20181223221242}
\end{eqnarray}
\textit{Step 2.}
 We prove that the variation $\mathcal{E}(\gamma)$ is finite when $\mu(d\xi)$ is infinite. For all $t>0$, let
\begin{eqnarray*}
a(t):=\log\mathbb{E}\exp
\left\{\frac{1}{2t}\int_0^{t}\int_0^{t}\gamma(B(s)-B(r))dsdr\right\},
\end{eqnarray*}
which is well defined because of the exponential integrability (\ref{11111}).
By Lemma \ref{C3.1}, we can obtain the subadditivity of $a(t)$, the definition of which is referred in \cite{R2}.
By the subadditivity and Lemma 1.3.4 in \cite{R2}, the following limit exists
\begin{eqnarray}\label{X32}
\Lambda:=\lim\limits_{t\rightarrow+\infty}\frac{1}{t}\log\mathbb{E}\exp
\left\{\frac{1}{2t}\int_0^{t}\int_0^{t}\gamma(B(s)-B(r))dsdr\right\}<+\infty.
\end{eqnarray}
We will use the same linearized method as the proof of lower bound in \textit{Step 1} and the $\mu_\varepsilon$, $\gamma_\varepsilon$ in (\ref{2019314224354}). For all $\varepsilon>0$, let the Hilbert space $\mathcal{H}_\varepsilon$ with $\mu_\varepsilon$, similarly to $\mathcal{H}$ in \textit{Step 1}, and $f^*_\varepsilon(x):=\int_{\mathrm{R}^d}e^{i\xi\cdot x}f(-\xi)\mu_\varepsilon(d\xi)$, we have
\begin{eqnarray}
&&\liminf\limits_{t\rightarrow+\infty}\frac{1}{t}\log\mathbb{E}\exp\left\{\frac{1}{2t}
\int_0^t\int_0^t\gamma_\varepsilon(B(s)-B(r))dsdr\right\}\nonumber\\
&&\ge\sup\limits_{f\in\mathcal{H}_\varepsilon}\left\{\liminf\limits_{t\rightarrow+\infty}\frac{1}{t}\log\mathbb{E}
\exp\left\{\int_0^{t}f^*_\varepsilon(B(s))ds\right\}-\frac{1}{2}\|f\|_{\mathcal{H}_\varepsilon}^2\right\}\nonumber\\
&&\ge\frac{1}{2}\sup\limits_{g\in\mathcal{F}_d}\left\{\int_{\mathrm{R}^d\times\mathrm{R}^d}\gamma_\varepsilon(x-y)g^2(x)g^2(y)dxdy-\int_{\mathrm{R}^d}|\nabla g(x)|^2dx\right\}.\label{X34}
\end{eqnarray}
For all $\varepsilon>0$, by Taylor expansion and Bochner's representation, we can get the inequality
\begin{eqnarray}\label{X33}
\mathbb{E}\exp\left\{\frac{1}{2t}\int_0^t\int_0^t\gamma_\varepsilon(B(s)-B(r))dsdr\right\}
\le\mathbb{E}\exp\left\{\frac{1}{2t}\int_0^t\int_0^t\gamma(B(s)-B(r))dsdr\right\}.
\end{eqnarray}
Then, by (\ref{X32})-(\ref{X33}) and monotone convergence theorem, we have
\begin{eqnarray*}
\mathcal{E}(\gamma)&=&\frac{1}{2}\sup\limits_{g\in\mathcal{F}_d}\left\{\int_{\mathrm{R}^d\times\mathrm{R}^d}\gamma(x-y)g^2(x)g^2(y)dxdy-\int_{\mathrm{R}^d}|\nabla g(x)|^2dx\right\}\\
&=&\frac{1}{2}\sup\limits_{\varepsilon>0}\sup\limits_{g\in\mathcal{F}_d}\left\{\int_{\mathrm{R}^d\times\mathrm{R}^d}\gamma_\varepsilon(x-y)g^2(x)g^2(y)dxdy-\int_{\mathrm{R}^d}|\nabla g(x)|^2dx\right\}\le\Lambda<+\infty.
\end{eqnarray*}

\textit{Step 3.} We complete the proof of the limit (\ref{X2}) when $\mu$ is infinite. For all $\varepsilon>0$,
by Taylor expansion and Bochner's representation, we can get the inequality
\begin{eqnarray}
&&\mathbb{E}\exp\Bigg\{\frac{1}{2Nt_N}\sum\limits_{j, k=1}^N\int_0^{t_N}\int_0^{t_N}\gamma_\varepsilon(B_j(s)-B_k(r))dsdr\Bigg\}\nonumber\\
&&\le\mathbb{E}\exp\Bigg\{\frac{1}{2Nt_N}\sum\limits_{j, k=1}^N\int_0^{t_N}\int_0^{t_N}\gamma(B_j(s)-B_k(r))dsdr\Bigg\}.\label{20181223223739}
\end{eqnarray}
Then, by (\ref{20181223223240}) and (\ref{20181223223739}), we get the lower bound as the followings
\begin{eqnarray*}
&&\liminf\limits_{N\rightarrow\infty}\frac{1}{Nt_N}\log\mathbb{E}\exp
\Bigg\{\frac{1}{2Nt_N}\sum\limits_{j, k=1}^N\int_0^{t_N}\int_0^{t_N}\gamma(B_j(s)-B_k(r))dsdr\Bigg\}\nonumber\\
&&\ge\sup\limits_{\varepsilon>0}\sup\limits_{g\in\mathcal{F}_d}\left\{\int_{\mathrm{R}^d\times\mathrm{R}^d}
\gamma_\varepsilon(x-y)g^2(x)g^2(y)dxdy-\frac{1}{2}\int_{\mathrm{R}^d}|\nabla g(x)|^2dx\right\}=\mathcal{E}(\gamma).
\end{eqnarray*}
For the upper bound, we make the decomposition  $\gamma=\gamma-\gamma_\varepsilon+\gamma_\varepsilon
:=\overline{\gamma}_{\varepsilon}
+\gamma_\varepsilon$, where $\overline{\gamma}_{\varepsilon}=\mathcal{F}(\overline{\mu}_{\varepsilon})$ in $\mathcal{S}'(\mathrm{R}^d)$ and $\overline{\mu}_{\varepsilon}(d\xi):=(1-\exp\left\{-\frac{1}{2}
\varepsilon|\xi|^2\right\})\mu(d\xi)$. For any parameters $ p,q>0$ and $ \frac{1}{p}+\frac{1}{q}=1$, by H\"{o}lder inequality, we have
\begin{eqnarray}
&&\mathbb{E}\exp\Bigg\{\frac{1}{2Nt_N}\sum\limits_{j, k=1}^N\int_0^{t_N}\int_0^{t_N}\gamma(B_j(s)-B_k(r))dsdr\Bigg\}\nonumber\\
&&\le\Bigg(\mathbb{E}\exp\Bigg\{\frac{p}{2Nt_N}\sum\limits_{j, k=1}^N\int_0^{t_N}\int_0^{t_N}\gamma_\varepsilon(B_j(s)-B_k(r))dsdr\Bigg\}
\Bigg)^{\frac{1}{p}}\label{20181223211634}\\
&&\cdot\Bigg(\mathbb{E}\exp\Bigg\{\frac{q}{2Nt_N}\sum\limits_{j, k=1}^N\int_0^{t_N}\int_0^{t_N}\overline{\gamma}_{\varepsilon}(B_j(s)-B_k(r))dsdr\Bigg\}
\Bigg)^{\frac{1}{q}}.\label{X35}
\end{eqnarray}
For the error (\ref{X35}), for any $q>0$, it suffices to prove
\begin{eqnarray}
\limsup\limits_{\varepsilon\rightarrow0}\limsup\limits_{N\rightarrow\infty}
\frac{1}{Nt_N}\log\mathbb{E}\exp
\Bigg\{\frac{q}{2Nt_N}\sum\limits_{j, k=1}^N\int_0^{t_N}\int_0^{t_N}\overline{\gamma}_{\varepsilon}(B_j(s)-B_k(r))
dsdr\Bigg\}\le 0.\label{X37}
\end{eqnarray}
Moreover, by Lemma \ref{C3.2} and the similar computations to (\ref{20181223221242}), we only need to prove
\begin{eqnarray}
&&\limsup\limits_{\varepsilon\rightarrow0}\limsup\limits_{N\rightarrow\infty}
\frac{1}{N}\log\mathbb{E}\exp
\Bigg\{\frac{q}{2N}\sum\limits_{j, k=1}^N\int_0^{1}\int_0^{1}\overline{\gamma}_{\varepsilon}(B_j(s)-B_k(r))
dsdr\Bigg\}\le 0.\label{2019320143120}
\end{eqnarray}
In fact, by the approximation in (\ref{2019314224354}), the similar computations to (\ref{201315101645}) and Jensen's inequality, we have
\begin{eqnarray}
\mathbb{E}\exp
\Bigg\{\frac{qN}{2}\int_{\mathrm{R}^{d}}\left|\frac{1}{N }\sum_{j=1}^N\int_{0}^{1}e^{i\xi\cdot B_j(s)}ds\right|^2\overline{\mu}_{\varepsilon}(d\xi)\Bigg\}\le\Bigg(\mathbb{E}\exp
\Bigg\{\frac{q}{2}\int_{\mathrm{R}^{d}}\left|\int_{0}^{1}e^{i\xi\cdot B(s)}ds\right|^2\overline{\mu}_{\varepsilon}(d\xi)\Bigg\}\Bigg)^{N}.\label{2019320162520}
\end{eqnarray}
By (\ref{2019320162520}), Taylor expansion and Fatou's Lemma, we find that (\ref{2019320143120}) comes from the following fact
\begin{eqnarray*}
\limsup\limits_{\varepsilon\rightarrow0}\mathbb{E}\exp\left\{\frac{q}{2}\int_{\mathrm{R}^{d}}\left|\int_{0}^{1}e^{i\xi\cdot B(s)}ds\right|^2\overline{\mu}_{\varepsilon}(d\xi)\right\}\le1.
\end{eqnarray*}
At last, by (\ref{20181223221242}) and (\ref{20181223211634})-(\ref{X37}), we have
\begin{eqnarray*}
&&\limsup\limits_{N\rightarrow\infty}\frac{1}{Nt_N}\log\mathbb{E}\exp\Bigg\{
\frac{1}{2Nt_N}\sum\limits_{j, k=1}^N\int_0^{t_N}\int_0^{t_N}\gamma(B_j(s)-B_k(r))dsdr\Bigg\}\\
&&\le\limsup\limits_{p\rightarrow1}\limsup\limits_{\varepsilon\rightarrow0}
\limsup\limits_{N\rightarrow\infty}\frac{1}{pNt_N}\log\mathbb{E}\exp\Bigg\{
\frac{p}{2Nt_N}\sum\limits_{j, k=1}^N\int_0^{t_N}\int_0^{t_N}\gamma_\varepsilon(B_j(s)-B_k(r))dsdr\Bigg\}\\
&&\le\lim\limits_{p\rightarrow1}
\frac{1}{2p}\sup\limits_{g\in\mathcal{F}_d}
\left\{p\int_{\mathrm{R}^d}\Bigg|\int_{\mathrm{R}^d}e^{i\xi\cdot x}g^2(x)dx\Bigg|^2\mu(d\xi)
-\int_{\mathrm{R}^d}|\nabla g(x)|^2dx\right\}=\mathcal{E}(\gamma).
\end{eqnarray*}

\section{An asymptotic result for asymptotically logarithmic function }\label{20181225200911}

To prove Theorem \ref{201946171643}, we may first consider the case of ``asymptotically logarithmic'' function, i.e. the function satisfying the following assumption, and prove a general result by Proposition \ref{C1}.
\begin{enumerate}
\item[ \noindent\textbf{(H)}]The $\gamma$ is point-wise defined in $\mathrm{R}^{d}\backslash\{0\}$ and bounded outside every neighborhood of $0$, and there exists some constant $C>0$ such that
\begin{eqnarray}
\qquad\qquad\gamma(x)\sim C\log\frac{1}{|x|}\quad\qquad(x\rightarrow0).\label{2019323162035}
\end{eqnarray}
\end{enumerate}
\begin{proposition}\label{2018109231801}
 Assume the condition (H), then for all $t>0$, it holds that
\begin{eqnarray}
\lim\limits_{N\rightarrow\infty}\frac{1}{N\sigma_N}
\log\mathbb{E}\exp\Bigg\{\frac{1}{2}\sum\limits_{ j,k=1}^N\int_0^t\int_0^t\gamma(B_j(s)-B_k(r))drds\Bigg\}=\frac{t}{2}.\label{201942104446}
\end{eqnarray}
Here, $\sigma_N$ is a function about positive integer $N$, which satisfies $\sigma_N>e$ and $\sigma_N^{-1}\log \sigma_N =\frac{2}{NCt}$ when $N$ is enough large.
\end{proposition}
\begin{remark}
 Proposition \ref{2018109231801} also includes the precise high moment asymptotic results for the PAM with massive free field and Bessel field. In fact, their covariances satisfy the condition (H),
because their covariances are Bessel potential $G_s(x)$ on $\mathrm{R}^d$ which approximates to the logarithmic potential $C_d\log\frac{1}{|x|}$ at $0$ with $C_d>0$ when $s=d$, by \cite{A23}. For the detailed definition of the two fields, see e.g. \cite{G22,R27,O26}.

\end{remark}

\noindent\textbf{The proof of Proposition \ref{2018109231801}:}
We will use the truncated power function $\phi_l(|x|):=(1-|x|)_+^l$ with $l>0$, which is positive definite on $\mathrm{R}^{d}$ when $l \ge\lfloor\frac{d}{2}\rfloor +1$ by \cite{H96}.
For any two parameters $\delta, M>0$, let $\gamma_\delta(x):=\phi_l(|\frac{x}{\delta}|)$ and $\gamma_M(x):=\phi_l(|\frac{x}{M}|)$.
By the condition (H), for any $\varepsilon>0$, there exists a $\delta>0$ such that for all $|x|<\delta$, it has
\begin{eqnarray}
(1-\varepsilon)C\log\frac{1}{|x|}
\le\gamma(x)\le (1+\varepsilon)C\log\frac{1}{|x|}.
\label{2018108204441}
\end{eqnarray}
To simplify notation, we let $C=1$ in (\ref{2018108204441}) and $(fg)(x):=f(x)g(x)$ for all $f,g$ on $\mathrm{R}^d$.

First, we prove the upper bound. Let $t_N=\sigma_N t$. For any parameters $ p,q>0$ and $\frac{1}{p}+\frac{1}{q}=1$, by Brownian scaling, $\gamma=\gamma\gamma_\delta+\gamma(1-\gamma_\delta)$ and H\"{o}lder inequality , we have
\begin{eqnarray}
&&\mathbb{E}\exp\Bigg\{\frac{1}{2}\sum\limits_{j, k=1}^N\int_0^{t}\int_0^{t}\gamma(B_j(s)-B_k(r))dsdr\Bigg\}\nonumber\\
&&\le\Bigg[\mathbb{E}\exp\Bigg\{\frac{p}{2\sigma_N^{2}}\sum\limits_{j, k=1}^N\int_0^{t_N}\int_0^{t_N}(\gamma
\gamma_\delta)(\sigma_N^{-\frac{1}{2}}(B_j(s)-B_k(r)))dsdr\Bigg\}\Bigg]^{\frac{1}{p}}\label{2018108205155}\\
&&\cdot\Bigg[\mathbb{E}\exp\Bigg\{\frac{q}{2\sigma_N^{2}}\sum\limits_{j, k=1}^N\int_0^{t_N}\int_0^{t_N}(\gamma
(1-\gamma_\delta))(\sigma_N^{-\frac{1}{2}}(B_j(s)-B_k(r)))dsdr\Bigg\}\Bigg]^{\frac{1}{q}}.\label{201818205208}
\end{eqnarray}
For the main term (\ref{2018108205155}),
by (\ref{2018108204441}) and $0\le\gamma_\delta\le1$, we have
\begin{eqnarray}
&&\mathbb{E}\exp\Bigg\{\frac{p}{2\sigma_N^{2}}\sum\limits_{j, k=1}^N\int_0^{t_N}\int_0^{t_N}(\gamma
\gamma_\delta)(\sigma_N^{-\frac{1}{2}}(B_j(s)-B_k(r)))dsdr\Bigg\}\qquad\qquad\qquad
\qquad\quad~~~\nonumber\\
&&\le\mathbb{E}\exp\Bigg\{\frac{(1+\varepsilon)p}{2\sigma_N ^{2}}\sum\limits_{j, k=1}^N
\int_0^{t_N}\int_0^{t_N} \log_+\frac{1}{|
(B_j(s)-B_k(r))|}dsdr+\frac{(1+\varepsilon)pt^2N^2
\log\sigma_N ^{\frac{1}{2}}}{2}\Bigg\}.\label{2018109153513}
\end{eqnarray}
Let $\sigma_N>e$ and $\frac{tN\log\sigma_N ^{\frac{1}{2}}}{\sigma_N}=1$ when $N$ is enough large, then $\sigma_N$ exists and satisfies $\frac{1}{\sigma_N}=o(\frac{1}{N})$ as $N\rightarrow\infty$. Then, by (\ref{2018109153513}), we have
\begin{eqnarray}
&&\limsup\limits_{N\rightarrow\infty}\frac{1}{pNt_N}\log\mathbb{E}\exp
\Bigg\{\frac{p}{2\sigma_N^{2}}\sum\limits_{j, k=1}^N\int_0^{t_N}\int_0^{t_N}(\gamma
\gamma_\delta)(\sigma_N^{-\frac{1}{2}}(B_j(s)-B_k(r)))dsdr\Bigg\}\nonumber\\
&&\le\limsup\limits_{N\rightarrow\infty}\frac{1}{pNt_N}\log\mathbb{E}
\exp\Bigg\{\frac{(1+\varepsilon)p}{2\sigma_N ^{2}}\sum\limits_{j, k=1}^N
\int_0^{t_N}\int_0^{t_N} \log_+\frac{1}{|
(B_j(s)-B_k(r))|}dsdr\Bigg\}+\frac{(1+\varepsilon)}{2}.\label{2018119184908}
\end{eqnarray}
We claim that for all $p>0$, it holds that
\begin{eqnarray}
\limsup\limits_{N\rightarrow\infty}\frac{1}{pNt_N}\log\mathbb{E}
\exp\Bigg\{\frac{(1+\varepsilon)p}{2\sigma_N ^{2}}\sum\limits_{j, k=1}^N
\int_0^{t_N}\int_0^{t_N} \log_+\frac{1}{|
(B_j(s)-B_k(r))|}dsdr\Bigg\}\le0.\label{2018109131601}
\end{eqnarray}
Indeed, to use Proposition \ref{C1}, we only need to find a positive definite function satisfying the Dalang's condition instead of $\log_+\frac{1}{|x|}$ in (\ref{2018109131601}), because the exact kernel $\log_+\frac{1}{|x|}$ is positive definite if and only if $d\le3$ in \cite{G29}. We may notice that there is some $\alpha\in(0,1)$ and $C>0$ such that
\begin{eqnarray}
\log_+\frac{1}{|x|}\le|x|^{-\alpha}+C.\label{20181019122223}
\end{eqnarray}
Furthermore, for all $\varepsilon_1>0$, it holds that $\frac{1}{\sigma_N ^{2}}=\frac{2}{\log\sigma_N }\frac{1}{Nt_N}\le\frac{ \varepsilon_1}{Nt_N}$ when $N$ is enough large.  By (\ref{20181019122223}) and Proposition \ref{C1}, we have
\begin{eqnarray*}
&&\limsup\limits_{N\rightarrow\infty}\frac{1}{pNt_N}\log\mathbb{E}
\exp\Bigg\{\frac{(1+\varepsilon)p}{2\sigma_N ^{2}}\sum\limits_{j, k=1}^N
\int_0^{t_N}\int_0^{t_N} \log_+\frac{1}{|
(B_j(s)-B_k(r))|}dsdr\Bigg\}\\
&&\le\limsup\limits_{\varepsilon_1\rightarrow0}\limsup\limits_{N\rightarrow\infty}
\frac{1}{pNt_N}\log\mathbb{E}
\exp\Bigg\{\frac{(1+\varepsilon)\varepsilon_1p}{2Nt_N}\sum\limits_{j, k=1}^N
\int_0^{t_N}\int_0^{t_N} (\left|B_j(s)-B_k(r)\right|^{-\alpha}+C)dsdr\Bigg\}\\
&&\le\limsup\limits_{\varepsilon_1\rightarrow0}
\frac{1}{2p}\sup\limits_{g\in\mathcal{F}_d}\Bigg\{(1+\varepsilon)\varepsilon_1
p\int_{\mathrm{R}^d\times\mathrm{R}^d}(|x-y|^{-\alpha}+C)
g^2(x)g^2(y)dxdy-\int_{\mathrm{R}^d}|\nabla g(x)|^2dx\Bigg\}=0.
\end{eqnarray*}
For the error term (\ref{201818205208}),
we only need to show that for all $q>0$, it holds that
\begin{eqnarray}
\limsup\limits_{N\rightarrow\infty}\frac{1}{qNt_N}\log\mathbb{E}
\exp\Bigg\{\frac{q}{2\sigma_N^{2}}\sum\limits_{j, k=1}^N\int_0^{t_N}\int_0^{t_N}(\left|\gamma
\right|(1-\gamma_\delta))(\sigma_N^{-\frac{1}{2}}(B_j(s)-B_k(r)))dsdr\Bigg\}
\le0.\label{2018109143109}
\end{eqnarray}
Indeed,
 for any $\varepsilon_2>0$, there exists a $\delta_1\in(0,\delta)$ such that $1-\gamma_\delta(x)\le\varepsilon_2$ when $|x|<\delta_1$. Since $\gamma$ is bounded outside every neighborhood of $0$, we can define $C:=\sup\limits_{|x|\ge\delta_1}|\gamma|<+\infty$. By (\ref{2018108204441}), we have
\begin{eqnarray}
|\gamma|(1-\gamma_\delta)
\le
(1+\varepsilon)\varepsilon_2\log\frac{1}{|x|}1_{|x|\le\delta_1}
+C.\label{2018109151444}
\end{eqnarray}
By (\ref{2018109151444}) and the similar computations to (\ref{2018109153513}), we obtain
\begin{eqnarray}
\mathbb{E}
\exp\Bigg\{\frac{q}{2\sigma_N^{2}}\sum\limits_{j, k=1}^N\int_0^{t_N}\int_0^{t_N}(\left|\gamma
\right|(1-\gamma_\delta))(\sigma_N^{-\frac{1}{2}}(B_j(s)-B_k(r)))dsdr\Bigg\}
\qquad\qquad\qquad\qquad\qquad\qquad\qquad~~\nonumber\\
\le\mathbb{E}
\exp\Bigg\{\frac{q(1+\varepsilon)\varepsilon_2}{2\sigma_N^{2}}\sum\limits_{j, k=1}^N\int_0^{t_N}\int_0^{t_N}\log_+\frac{1}{|B_j(s)-B_k(r)|}dsdr+\frac{1}{2}q(1+\varepsilon)\varepsilon_2t^2N^2
\log\sigma_N ^{\frac{1}{2}}+ \frac{1}{2}Cqt^2N^2\Bigg\}.\label{2018119221405}
\end{eqnarray}
Moreover, by the same computations as (\ref{2018109153513})-(\ref{2018109131601}) and (\ref{2018119221405}), we have
\begin{eqnarray*}
&&\limsup\limits_{N\rightarrow\infty}\frac{1}{qNt_N}\log\mathbb{E}
\exp\Bigg\{\frac{q}{2\sigma_N^{2}}\sum\limits_{j, k=1}^N\int_0^{t_N}\int_0^{t_N}(\left|\gamma
\right|(1-\gamma_\delta))(\sigma_N^{-\frac{1}{2}}(B_j(s)-B_k(r)))dsdr\Bigg\}\nonumber\\
&&\le\limsup\limits_{\varepsilon_2\rightarrow0}\frac{(1+\varepsilon)
\varepsilon_2}{2}=0.
\end{eqnarray*}
Then, by (\ref{2018108205155}), (\ref{201818205208}), (\ref{2018119184908}), (\ref{2018109131601}) and (\ref{2018109143109}), we get the upper bound
\begin{eqnarray*}
\limsup\limits_{N\rightarrow\infty}\frac{1}{Nt_N}\log\mathbb{E}\exp\Bigg\{
\frac{1}{2}\sum\limits_{j, k=1}^N\int_0^{t}\int_0^{t}\gamma(B_j(s)-B_k(r))dsdr\Bigg\}\le\limsup\limits_{\varepsilon\rightarrow0}\frac{(1+\varepsilon)}{2}=\frac{1}{2}.
\end{eqnarray*}
Second, we prove the lower bound. For any parameter $ p>1$, by the reverse H\"{o}lder inequality, we have
\begin{eqnarray}
&&\mathbb{E}\exp\Bigg\{\frac{1}{2}\sum\limits_{j, k=1}^N\int_0^{t}\int_0^{t}\gamma(B_j(s)-B_k(r))dsdr\Bigg\}\nonumber\\
&&\ge\Bigg[\mathbb{E}\exp\Bigg\{\frac{p}{2\sigma_N^{2}}\sum\limits_{j, k=1}^N\int_0^{t_N}\int_0^{t_N}(\gamma
\gamma_\delta)(\sigma_N^{-\frac{1}{2}}(B_j(s)-B_k(r)))
dsdr\Bigg\}\Bigg]^{\frac{1}{p}}\label{2018109173452}\\
&&\cdot\Bigg[\mathbb{E}\exp\Bigg\{\frac{-1}{2(p-1)\sigma_N^{2}}\sum\limits_{j, k=1}^N\int_0^{t_N}\int_0^{t_N}(\gamma
(1-\gamma_\delta))(\sigma_N^{-\frac{1}{2}}(B_j(s)-B_k(r)))
dsdr\Bigg\}\Bigg]^{(1-p)}.\label{2018109173505}
\end{eqnarray}
For any $M>1$, there exists a $N_0>0$ such that for all $N>N_0$, it has $M\le\delta\sigma_N^{\frac{1}{2}}$.
When $N>N_0$, by (\ref{2018108204441}) and $0\le\gamma_M\le1$, we get the following evaluation for (\ref{2018109173452})
\begin{eqnarray}
&&\mathbb{E}\exp\Bigg\{\frac{p}{2\sigma_N^{2}}\sum\limits_{j, k=1}^N\int_0^{t_N}\int_0^{t_N}(\gamma
\gamma_\delta)(\sigma_N^{-\frac{1}{2}}(B_j(s)-B_k(r)))
dsdr\Bigg\}\nonumber\\
&&\ge
\mathbb{E}\exp\Bigg\{\frac{(1-\varepsilon)p}{2\sigma_N ^{2}}\log\sigma_N ^{\frac{1}{2}}\sum\limits_{j, k=1}^N
\int_0^{t_N}\int_0^{t_N}\gamma_M(B_j(s)-B_k(r))dsdr\nonumber\\
&&+\frac{(1-\varepsilon)p}{2\sigma_N ^{2}}\sum\limits_{j, k=1}^N
\int_0^{t_N}\int_0^{t_N} \log\frac{1}{|
(B_j(s)-B_k(r))|}\gamma_M(B_j(s)-B_k(r))dsdr\Bigg\}\nonumber\\
&&\ge\mathbb{E}\exp\Bigg\{\frac{(1-\varepsilon)p}{2\sigma_N ^{2}}\log\sigma_N ^{\frac{1}{2}}\sum\limits_{j, k=1}^N
\int_0^{t_N}\int_0^{t_N}\gamma_M(B_j(s)-B_k(r))dsdr-p\log M\frac{t^2N^2}{2}
\Bigg\}.\label{2018109202532}
\end{eqnarray}
Moreover, by (\ref{2018109202532}), Proposition \ref{C1} and $\frac{\log\sigma_N^{\frac{1}{2}}}{\sigma_N^2}=\frac{1}{Nt_N}$ (i.e. $\frac{t^2N}{t_N}=\frac{2}{\log\sigma_N}$) when $N$ is enough large, we have
\begin{eqnarray}
&&\liminf\limits_{N\rightarrow\infty}\frac{1}{pNt_N}\log\mathbb{E}
\exp\Bigg\{\frac{p}{2\sigma_N^{2}}\sum\limits_{j, k=1}^N\int_0^{t_N}\int_0^{t_N}(\gamma
\gamma_\delta)(\sigma_N^{-\frac{1}{2}}(B_j(s)-B_k(r)))
dsdr\Bigg\}\nonumber\\
&&\ge\liminf\limits_{N\rightarrow\infty}\frac{1}{pNt_N}\log\mathbb{E}\exp\Bigg
\{\frac{(1-\varepsilon)p}{2Nt_N}\sum\limits_{j, k=1}^N
\int_0^{t_N}\int_0^{t_N}\gamma_M(B_j(s)-B_k(r))dsdr\Bigg\}\nonumber\\
&&\ge\frac{1}{2p}\sup\limits_{M>1}\sup\limits_{g\in\mathcal{F}_d}\Bigg\{(1-\varepsilon)p
\int_{\mathrm{R}^d\times\mathrm{R}^d}\gamma_M(x-y)
g^2(x)g^2(y)dxdy-\int_{\mathrm{R}^d}|\nabla g(x)|^2dx\Bigg\}=\frac{1-\varepsilon}{2}.\label{2018119232030}
\end{eqnarray}
Here, the last step is due to $\gamma_M\uparrow1$ as $M\rightarrow+\infty$.
For the error term (\ref{2018109173505}), we only need to show that for all $p>1$, it holds that
\begin{eqnarray}
\liminf\limits_{N\rightarrow\infty}
\frac{-(p-1)}{Nt_N}\log\mathbb{E}\exp\Bigg\{\frac{-1}{2(p-1)\sigma_N^{2}}\sum\limits_{j, k=1}^N\int_0^{t_N}\int_0^{t_N}(\gamma
(1-\gamma_\delta))
(\sigma_N^{-\frac{1}{2}}(B_j(s)-B_k(r)))dsdr\Bigg\}\ge0.\label{2018109225613}
\end{eqnarray}
Indeed, for all $p>1$, by (\ref{2018109143109}), we have
\begin{eqnarray*}
&&\limsup\limits_{N\rightarrow\infty}
\frac{(p-1)}{Nt_N}\log\mathbb{E}\exp\Bigg\{\frac{-1}{2(p-1)\sigma_N^{2}}
\sum\limits_{j, k=1}^N\int_0^{t_N}\int_0^{t_N}(\gamma
(1-\gamma_\delta))(\sigma_N^{-\frac{1}{2}}(B_j(s)
-B_k(r)))dsdr\Bigg\}\nonumber\\
&&\le\limsup\limits_{N\rightarrow\infty}\frac{(p-1)}{Nt_N}\log\mathbb{E}
\exp\Bigg\{\frac{1}{2(p-1)\sigma_N^{2}}\sum\limits_{j, k=1}^N\int_0^{t_N}\int_0^{t_N}(\left|\gamma
\right|(1-\gamma_\delta))(\sigma_N^{-\frac{1}{2}}(B_j(s)
-B_k(r)))dsdr\Bigg\}=0.
\end{eqnarray*}
At last, by (\ref{2018109173452}), (\ref{2018109173505}), (\ref{2018119232030}) and (\ref{2018109225613}), we get the lower bound
\begin{eqnarray*}
~~~~~~\liminf\limits_{N\rightarrow\infty}\frac{1}{Nt_N}\log\mathbb{E}
\exp\Bigg\{\frac{1}{2}\sum\limits_{j, k=1}^N\int_0^{t}\int_0^{t}\gamma(B_j(s)-B_k(r))dsdr\Bigg\}\ge\liminf
\limits_{\varepsilon\rightarrow0}\frac{(1-\varepsilon)}{2}=\frac{1}{2}.
\end{eqnarray*}
\section{The proof of Theorem \ref{201946171643}}\label{201947120052}

Because the kernel $k(x,y)$ may be non-stationary, we need to find a function satisfying the condition (H) instead of the $k(x,y)$ in the moment representation (\ref{20181226223826}).
By (\ref{20181018215909}), there exists a constant $C>0$ such that $\log_+\frac{1}{|x-y|}-C\le k(x,y)\le \log_+\frac{1}{|x-y|}+C$ on $\mathrm{R}^{d}\times\mathrm{R}^{d}$. Furthermore,
by (\ref{20181217092348}) and
$\frac{1}{N\lambda_N}=o(\frac{1}{N^2})$ as $N\rightarrow\infty$ in Remark \ref{2019412145239}, we notice that when $N\rightarrow\infty$, it has
\begin{eqnarray}
\frac{1}{N\lambda_N}\log\mathbb{E}u^N(t,x)
\sim\frac{t}{N\sigma_N}\log\mathbb{E}\exp\Bigg\{\frac{1}{2}\sum\limits_{ j,k=1}^N\int_0^t\int_0^t\log_+\frac{1}{|B_j(s)-B_k(r)|}drds\Bigg\}.\label{20181216150115}
\end{eqnarray}
In fact, to simplify the notation, we use the $\lambda_N$ instead of the $\sigma_N$ in Proposition \ref{2018109231801}. By $\frac{\log \lambda_N}{\lambda_N} =\frac{2}{N}$ and $\frac{\log \sigma_N}{\sigma_N} =\frac{2}{Nt}$ and $\lambda_N,\sigma_N<N^2$ when $N$ is enough large, we have $\lim\limits_{N\rightarrow\infty}\frac{\sigma_N}{\lambda_N}
=\lim\limits_{N\rightarrow\infty}\frac{t\log(Nt/2\log\sigma_N)}
{\log(N/2\log\lambda_N)}
=t$.
By (\ref{20181216150115}), Theorem \ref{201946171643} is a direct result of Proposition \ref{2018109231801}.

\section{Acknowledgments}

The author would like to thank Professor Xia Chen for his help during the completion of this paper. The manuscript is partially supported by NSFC Grant No. 11871244.


\begin{thebibliography}{99}{\small
\bibitem{F1}
X. Chen, T. V. Phan, Free energy in a mean field of Brownian particles, Disrete Cont. Dyn. S. 39(2) (2019) 747-769.
\bibitem{R2}
X. Chen, Random Walk Intersections: Large Deviations and Related Topics, Math. Surveys Monogr., vol. 157, Amer.
Math. Soc., Providence, RI, 2010.
\bibitem{D5}
 M.D. Donsker, S.R.S. Varadhan, Asymptotic evaluation of certain Markov process expectations for large time III, Comm. Pure Appl. Math. 29 (1976) 389-461.
\bibitem{D6}
A. Dembo and O. Zeitouni, Large Deviations Techniques and Applications (2nd ed.), Springer, New York, 1998.
\bibitem{S3}
 R. A. Carmona, S. A. Molchanov, Stationary parabolic Anderson model and  intermittency, Probab. Theory Related Fields 102 (1995) 433-453.
\bibitem{M7}
J. G\"{a}rtner, W. K\"{o}nig, Moment asymptotics for the continuous parabolic Anderson model, Ann. Appl. Probab. 10 (2000) 192-217.
\bibitem{Q20}
X. Chen, Quenched asymptotics for Brownian motion in generalized Gaussian potential, Ann. Probab. 42 (2014) 576-622.
\bibitem{O9}
D. Conus,  M. Joseph, D. Khoshnevisan, On the chaotic character of the  stochastic heat equation, before the onset of intermitttency, Ann. Probab. 41 (2013) 2225-2260.
\bibitem{M10}
X. Chen, Moment asymptotics for parabolic Anderson equation with fractional time-space noise: in Skorokhod regime, Annales de l'Institut Henri Poincar\'{e} 53 (2017) 819-841.
\bibitem{S11}
X. Chen, Spatial asymptotics for the parabolic Anderson models with generalized time-space
Gaussian noise, Ann. Probab. 44 (2016) 1535-1598.
\bibitem{S14}
Y. Hu, J. Huang, D. Nualart, S. Tindel, Stochastic heat equations with general multiplicative Gaussian noises: H\"{o}lder continuity and intermittency, Electron. J. Probab. 20 (2015), no. 55, 50.
\bibitem{F15}
Y. Hu, D. Nualart, J. Song, Feynman-Kac formula for heat equation driven by fractional white noise, Ann. Probab. 39 (2011) 291-326.
\bibitem{P16}
X. Chen, Parabolic Anderson model with rough or critical Gaussian noise, Annales de l'Institut Henri Poincar\'{e} (accepted), \url{http://www.math.utk.edu/~xchen/rough-6.pdf}.
\bibitem{W18}
W. K\"{o}nig, The Parabolic Anderson Model: Random Walk in Random Potential, Birkh\"{a}user, 2016.
\bibitem{S19}
X. Chen, Y. Hu, D. Nualart, S. Tindel, Spatial asymptotics for the parabolic Anderson model driven by a Gaussian rough noise, Electron. J. Probab. 22 (2017), no. 65, 38.
\bibitem{P21}
H. Li, X. Chen, Precise moment asymptotics for the stochastic heat equation of a time-derivative Gaussian noise, Acta. Math. Sci. (English ser) (accepted), \url{http://www.math.utk.edu/~xchen/variance-1.pdf}.
\bibitem{G22}
 S. Sheffield, Gaussian free fields for mathematicians, Probab. Theory Related Fields, 139(3-4) (2007) 521-541.
\bibitem{A23}
 N. Aronszajn, K. T. Smith, Theory of Bessel potentials I, Ann. Inst. Fourier. 11 (1961) 385-475.
\bibitem{S25}
B. Duplantier,  R. Rhodes, S. Sheffield, V. Vargas, Critical Gaussian multiplicative chaos: convergence of the derivative martingale, Ann. Probab. 42(5) (2012) 1769-1808.
\bibitem{M25}
T. Madaule, Maximum of a log-correlated Gaussian field, Annales de l'Institut Henri Poincar\'{e} 51(4) (2015) 1369-1431.
\bibitem{R27}
B. Duplantier, R. Rhodes, S. Sheffield, V. Vargas, Renormalization of critical Gaussian multiplicative chaos and KPZ relation, Commun. Math. Phys. 330(1) (2014) 283-330.
\bibitem{L31}
R. Rhodes, V. Vargas, Lectures on Gaussian Multiplicative Chaos, \url{http://www.newton.ac.uk/files/seminar/20150119100011001-297514.pdf}.
\bibitem{O26}
L. Pitt, R. Robeva, On the sharp Markov property for Gaussian random fields and spectral synthesis in spaces of Bessel potentials, Ann. Probab. 31(3) (2003) 1338-1376.
\bibitem{G29}
R. Robert, V. Vargas, Gaussian multiplicative chaos revisited, Ann. Probab.  38(2) (2010) 605-631.
\bibitem{Q35}
 P. Chakraborty, X. Chen, B. Gao, S. Tindel, Quenched asymptotics for a 1-D stochastic heat equation driven by a rough spatial noise (preprint), \url{https://arxiv.org/abs/1810.04212}.
\bibitem{H96}
H. Wendland, Scattered Data Approximation (Cambridge Monographs on Applied and Computational Mathematics) 1st Edition, Cambridge University Press, 2004.


}
\end{thebibliography}
\end{document}